\documentclass[12pt]{article}

\usepackage{comment}
\usepackage{amsmath, amssymb, latexsym, amsthm}
\usepackage{color}
\usepackage{graphicx}
\usepackage{stmaryrd,amsbsy,enumerate}
\usepackage{hyperref}
\usepackage{enumitem}

\newtheorem{theorem}{Theorem} 
\newtheorem{theorem*}{Theorem}

\newtheorem{prop}[theorem]{Proposition}
\newtheorem{problem}{Problem}

\let\oldlfloor\lfloor
\let\oldrfloor\rfloor
\let\oldlceil\lceil
\let\oldrceil\rceil
\renewcommand{\lfloor}{\left\oldlfloor}
\renewcommand{\rfloor}{\right\oldrfloor}
\renewcommand{\lceil}{\left\oldlceil}
\renewcommand{\rceil}{\right\oldrceil}

\newcommand{\eps}{\varepsilon}

\newcommand{\cG}{\mathcal{G}}

\DeclareMathOperator\mydeg{deg}
\DeclareMathOperator\polylog{polylog}
\newcommand{\hrefemail}[1]{\href{#1}{#1}}

\begin{document}

\title{Subgraphs with a positive minimum semidegree in digraphs with large outdegree}
\author{
Andrzej Grzesik~\thanks{
Faculty of Mathematics and Computer Science, Jagiellonian University, ul.~Prof.~St.~\L ojasiewicza 6, 30-348 Krak\'{o}w, Poland.
E-mail: \hrefemail{Andrzej.Grzesik@uj.edu.pl}.
This author was supported by the National Science Centre grant 2021/42/E/ST1/00193.}
\and
Vojtěch Rödl~\thanks{
Department of Mathematics, Emory University, Atlanta, USA.
E-mail: \hrefemail{vrodl@emory.edu}.
This author was supported by the National Science Foundation (DMS~2300347).}
\and
Jan Volec~\thanks{
Department of Theoretical Computer Science, Faculty of Information Technology,
Czech Technical University in Prague, Th\'akurova 9, Prague, 160 00, Czech Republic.
E-mail: \hrefemail{jan@ucw.cz}.
This author was supported by the grant 23-06815M of the Grant Agency of the Czech Republic.}
}

\date{}
\maketitle

\begin{abstract}
We prove that every $n$-vertex directed graph $G$ with the minimum outdegree $\delta^+(G) = d$ contains a subgraph $H$ satisfying 
\[
\min\left\{\delta^+(H), \delta^-(H) \right\} \ge \frac{d(d+1)}{2n}
\,.\]
We also show that if $d = o(n)$ then this bound is asymptotically best possible.
\end{abstract}

\section{Introduction}\label{sec:intro}
 
A well known result of Erd\H{o}s~\cite{bib:E6465} states that every undirected graph with $n$~vertices and $m$ edges contains a subgraph with minimum degree at~least~$\frac mn$.
In order to establish a meaningful analogue of this statement in the setting of directed graphs, one needs to impose some additional assumptions;
indeed every subgraph of a transitive tournament contains a vertex with zero indegree (and similarly a vertex with zero outdegree).
Huang et al.~\cite{bib:HMSSY13} proved that if a directed graph with $n$ vertices and $m$ arcs is Eulerian, i.e., if every vertex has its indegree equal to its outdegree,
then one can always find an~Eulerian subgraph with the minimum indegree at least $\frac{m^2}{24n^3}$.

In this paper, we focus on directed graphs where every vertex has a~large outdegree.
It is easy to observe that if every vertex in a finite digraph has non-zero outdegree, then it has to contain a directed cycle.
An immediate corollary of this observation is that such a digraph contains a subgraph where both the minimum indegree and the minimum outdegree are non-zero.
Our aim is to study a quantitative strengthening of the aforementioned corollary.

\subsection{Notation}\label{sec:notation}

Given a finite vertex-set $V$, a \emph{directed graph}, or \emph{digraph} for short, is a pair $(V,A)$, where $A \subseteq V\times V$.
The elements of $A$ are called \emph{arcs}, and for~a~given arc $e=(u,w)$, where $u,w\in V$, we say that $e$ is \emph{oriented} from $u$ to $w$.
We note that all the digraphs in this paper are \emph{loopless}, i.e., \hbox{$A \subseteq V \times V \setminus \{(v,v): v\in V\}$},
however, pairs of different vertices are allowed to be connected by two arcs with the opposite directions.
For a systematic study of digraphs, we refer the reader to the monograph of Bang-Jensen and Gutin~\cite{bib:B-JG}.

Given a digraph $G=(V,A)$ and a vertex $v \in V$, we define the \emph{outneighborhood of $v$} and \emph{inneighborhood of $v$} to be \[
N_G^+(v) := \{ w \in V: (v,w) \in A\} \quad \mbox{and}\quad N_G^-(v) := \{ u \in V: (u,v) \in A\}
\,,  \]
respectively.
We set $\mydeg^+_G(v) := |N^+_G(v)|$ and $\mydeg^-_G(v):=|N_G^-(v)|$, and call the respective quantities the \emph{outdegree of $v$} and the \emph{indegree of $v$}.
We drop the subscript in cases when the digraph $G$ is clear from the context.

For a digraph $G$ on the vertex-set $V$, we define the~\emph{minimum outdegree of $G$}, which we denote by $\delta^+(G)$, to be the smallest outedgree of a vertex in $G$.
In other words, $\delta^+(G) := \min_{v \in V} \mydeg^+_G(v)$.
Analogously, we define the~\emph{minimum indegree of $G$} as $\delta^-(G) := \min_{v \in V} \mydeg^-_G(v)$.
Finally, the~\emph{minimum semidegree of $G$}, which we denote by $\delta^\pm(G)$, is the minimum of $\delta^+(G)$ and $\delta^-(G)$.

\subsection{Our results}

Using the notation from Section~\ref{sec:notation}, the observation from the beginning can be reformulated as follows:
every digraph $G$ with non-zero $\delta^+(G)$ contains a subgraph $H$ with $\delta^\pm(H)\ge1$.
It is natural to ask whether a stronger lower bound on $\delta^+(G)$ guarantees a stronger bound on $\delta^\pm(H)$.
The following example shows that even an assumption $\delta^+(G) = \sqrt{2n}-2$ is not sufficient to~guarantee a subgraph $H$ in $G$ with $\delta^\pm(H) \ge 2$:
\begin{enumerate}
\item Partition the set $[n]$ into two sets $A:=\left[\sqrt{2n}\right]$ and $B:=[n]\setminus A$.
\item Orient all the pairs $i<j$ in $A$ from $j$ to $i$.
\item For every $i \in A$, select $\sqrt{2n}-i-1$ outneighbors of $i$ in $B$ such a way that every $v \in B$ has at most one inneighbor in $A$;
note that this can be done since $\sum_{i=1}^{\sqrt{2n}} \sqrt{2n}-i-1 = n-\sqrt{4.5n} < |B|$.
\item Orient all the remaining pairs between $A$ and $B$ from $B$ to $A$.
\end{enumerate}
A moment of thought reveals that every vertex has at least $\sqrt{2n}-2$ outneighbors, yet every subgraph $H$ satisfies $\delta^\pm(H) \le 1$;
see Section~\ref{sec:construction} for details.
The~main result of this paper is that a lower bound $\delta^+(G) \ge \sqrt{2n}$ already yields a subgraph $H$ with $\delta^\pm(H) > 1$.
\begin{theorem}\label{thm:main}
Every loopless digraph $G$ with $n$ vertices and outdegree $\delta^+(G) \ge d$ has a subgraph $H$ with $\delta^\pm(H) \ge \frac{d(d+1)}{2n}$.
\end{theorem}

The rest of the paper is organized as follows.
We prove Theorem~\ref{thm:main} in Section~\ref{sec:proof}.
In Section~\ref{sec:construction}, we generalize the construction from the previous paragraph and show that Theorem~\ref{thm:main} is  for $d = o(n)$,
up to lower order terms, best possible.
We conclude by Section~\ref{sec:outro}, where we discuss the situation for the dense case $d = \alpha n$, and state some related open problems.

\section{Proof of Theorem~\ref{thm:main}}\label{sec:proof} 

Fix an $n$-vertex digraph $G$.
Without loss of generality, we will assume that every vertex of $G$ has the outdegree equal to $d$, as otherwise we pass to an appropriate spanning subgraph of $G$.
Now consider the following greedy procedure of one-by-one removing vertices from $G$:
\begin{enumerate}
\item Let $G_n := G$.
\item For $i = n,n-1,\ldots,2,1$, do:
\begin{enumerate}[leftmargin=\parindent]
\item Let $v_i$ be a vertex in $G_i$ with $\min\{\deg_{G_i}^+(v_i),\deg_{G_i}^-(v_i)\} = \delta^\pm(G_i)$.
\item Set $G_{i-1} := G_i - v_i$.
\end{enumerate}
\end{enumerate}
We note that in case of a tie between the minimum indegree and the minimum outdegree during the $i$-th step,
i.e., when $\delta^+(G_i) = \delta^-(G_i)$, we choose as $v_i$ a vertex that has the smallest indegree.

For the rest of the proof, we will be considering the vertices of $G$ being ordered as $v_1 < v_2 < \ldots < v_n$.
In other words, it is the reversed order with respect to the removal during the greedy procedure.
For brevity, given a vertex $v \in V(G)$, we will be writing $N^+_L(v)$ and $N^+_R(v)$
to denote the set of the outneighbors of $v$ in $G$ that are in the order before $v$ and after $v$, respectively.
Analogously, we define $N^-_L$ and $N^-_R$. In other words,
\begin{align*}
N^+_L(v) := & \left\{ w \in N_G^+(v): w<v \right\}\,, & N^+_R(v) := & \left\{ w \in N_G^+(v): w>v \right\}\,, \\
N^-_L(v) := & \left\{ w \in N_G^-(v): w<v \right\}\,, & N^-_R(v) := & \left\{ w \in N_G^-(v): w>v \right\}\,.
\end{align*}

Let $V^- := \left\{v_i : \mydeg^-_{G_i}(v_i) = \delta^\pm(G_i) \right\}$, and $V^+ := V(G) \setminus V^-$.
In words, the set $V^-$ consists of all the vertices that the greedy procedure removed because of their small indegree,
and the vertices in $V^+$ are those that were removed because of their small outdegree.
Set $c := \max_{i\in[n]} \delta^\pm(G_i)$.
Clearly, it suffices to show that $c \ge \frac{d(d+1)}{2n}$.
We will do so by distinguishing two cases based on the size of $V^+$.

We first focus on the case $|V^+| \ge d$.
Let $X$ be the last $d$ vertices from $V^+$, i.e., the first $d$ vertices that the greedy procedure removed because of having the smallest outdegree in the corresponding subgraph.
Since every vertex  of~the~digraph $G$ has outdegree equal to $d$,
it holds that $d^2 = \sum_{x\in X} \mydeg_G^+(x)$.
On the other hand, every vertex $x \in X$ satisfies
\[ \mydeg^+_G(x) = |N^+_L(x)| + |N^+_R(x) \cap X| + |N^+_R(x) \setminus X| \,.\]
By the definition of $X \subseteq V^+$, every $x\in X$ also satisfies $|N^+_L(x)| \le c$.
Furthermore, we claim that
\begin{equation}\label{eq:caseA:inside}
\sum_{x\in X} |N^+_R(x)\cap X| \le \binom{d}2 
\,.
\end{equation}
Indeed, every pair of the vertices from $X$ can contribute to this sum by at most one.
Finally, the fact that $X$ consists of the rightmost vertices of $V^+$ implies that $\bigcup\limits_{x \in X} \left(N_R^+(x) \setminus X\right) \subseteq V^-$,
which together with a standard double-counting argument yields that
\begin{equation*}\label{eq:caseA:between}
\sum_{x\in X} |N^+_R(x) \setminus X| = \sum_{x\in X} |N^+_R(x) \cap V^-| = \sum_{w \in V^-} |N^-_L(w) \cap X| \le (n-d)c
\,.
\end{equation*}
Therefore, we conclude that
\[
d^2 = \sum_{x\in X} |N^+_L(x)| + |N^+_R(x) \cap X| + |N^+_R(x) \setminus X|
\le d \cdot c + \binom{d}{2} + (n-d) c
\,.\]
Rearranging the last inequality readily yields that $c \ge \frac{d(d+1)}{2n}$.

It remains to consider the case when $|V^+| < d$.
This time, let $Y^-$ be the~first $\left(d - |V^+|\right)$ vertices from $V^-$, and let $Y:=Y^- \cup V^+$.
Firstly, observe that all the vertices to the left of a vertex $y\in Y^-$ must lie inside the set $Y$.
In particular, $N_L^-(y) \subseteq Y$ and $N_L^+(y) \subseteq Y$ for any $y\in Y^-$.
Analogously as~in~\eqref{eq:caseA:inside}, it holds that
\begin{equation*}\label{eq:caseB:inside}
\sum_{y\in Y^-} |N^+_L(y)| + \sum_{y\in Y} |N^+_R(y) \cap V^+| \le \binom{d}2
\,,
\end{equation*}
since the left-hand side of the inequality counts every pair of the vertices from $Y$ at most once.
We also have that $|N^+_L(y)| \le c$ for every $y \in V^+$, and additionally, it holds that $|N^-_L(y)| \le c$ for every $y \in Y^-$.
Thus,
\begin{equation*}\label{eq:caseB:c_controlled}
\sum_{y \in V^+} |N^+_L(y)| + \sum_{y\in Y} |N^+_R(y) \cap Y^-|
\le
|V^+|\cdot c + \sum_{w\in Y^-} |N^-_L(w)|
\le  d\cdot c
\,.
\end{equation*}
Finally, a double-counting argument yields that
\begin{equation*}\label{eq:caseB:between}
\sum_{y \in Y} |N^+_R(y) \setminus Y| = \sum_{w \in V \setminus Y} |N^-_L(w) \cap Y| \le (n-d)c
\,.
\end{equation*}
Since \[
|N_R^+(y)| = \left|N_R^+(y) \cap V^+\right| + \left|N_R^+(y) \cap Y^-\right| + \left|N_R^+(y) \setminus Y\right| \mbox{ for all } y \in Y,
\] we have that \[
d^2 = \sum_{y \in Y} |N^+_R(y)| + \sum_{y\in Y^-} |N_L^+(y)| + \sum_{y \in V^+} |N_L^+(y)| \le \binom{d}2 + d\cdot c +  (n-d)c
.\]
As in the previous case, rearranging $d^2 \le n \cdot c + \binom d2$ yields that $c \ge \frac{d(d+1)}{2n}$.

\section{\large Tournaments with all subgraphs having small~\texorpdfstring{$\delta^\pm$}{semidegree}}\label{sec:construction}
In this section, we prove that Theorem~\ref{thm:main} is essentially best possible, even if we would consider only \emph{tournaments}, i.e., orientations of the edge-set of complete undirected graphs.

\begin{prop}\label{prop:construction}
Fix two positive integers $k$ and $\ell$, and let $d:=2k\ell-1$ and $n_0:=(k+1)d+\ell$.
For every integer $n \ge n_0$ there exists a tournament $T$ with $n$ vertices and $\delta^+(T)=d$ such that
\[\max_{H\subseteq T} \delta^\pm(H) \le \ell \le \left(1 + \frac{k+1}{k^2}\right) \cdot \frac{d(d+1)}{2n}.\]
\end{prop}

\begin{figure}[t]
\centering \includegraphics[scale=0.9]{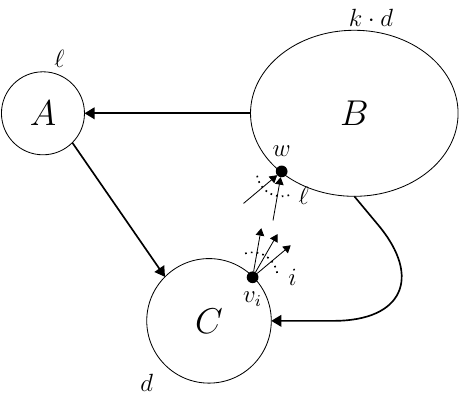}
\caption{The tournament constructed in Proposition~\ref{prop:construction}.}
\end{figure}

\begin{proof}
Fix the parameters $k$ and $\ell$.
Firstly, let us consider the case $n = n_0$, and let $V$ be any vertex-set of cardinality $n$.
We partition $V$ into three parts $A$, $B$ and $C$ of respective sizes $\ell$, $k\cdot d$ and $d$.
Let us now describe the orientations of pairs of vertices from $V$.
Firstly,
\begin{enumerate}
\item each of the sets $A$, $B$, $C$ induces a transitive orientation,
\item all the pairs between $A$ and $B$ are oriented from $B$ to $A$, and
\item all the pairs between $A$ and $C$ are oriented from $A$ to $C$.
\end{enumerate}
Clearly, $\deg^+_T(v) \ge |C| = d$ for every vertex $v \in A$.
Now, in order to have $\mydeg_T^+(v) = d$ for every $v \in C$,
we need to orient a total of $\binom{d+1}2$ pairs from $C$ to $B$,
where the $i$-th vertex of $C$ with respect to the transitive order needs to have exactly $i$ outneighbors in $B$.
Our aim is to do this in such a way that every vertex of $B$ will have exactly $\frac{d+1}{2k} = \ell$ incoming arcs from~$C$.
We claim that this is always possible. 
Indeed, order $B$ arbitrarily and let $S$ be a sequence of length $\ell\cdot|B|$ formed by $\ell$ copies of the ordered $B$ placed after each other.
Since $\ell |B| = \binom{d+1}2$ and $|B| \ge d$, we can go over the vertices in $C$ one by one and orient arcs from the $i$-th vertex to the next $i$ members of $S$.
This establishes the desired $\mydeg_T^+(v) = d$ for every $v \in C$ as well as $|N^-(w) \cap C| = \ell$ for every $w \in B$.
We orient all the remaining pairs between $B$ and $C$ from $B$ to~$C$.
In particular, every vertex in $B$ has now exactly $d-\ell + |A| = d$ outneighbors in $C \cup A$.

In the previous paragraph we have shown that $\delta^+(T) = d$,
thus it remains to find an upper bound on $\delta^\pm(H)$ for any given $H \subseteq T$.
Without loss of generality, we will assume that $H$ is induced.
If $V(H)$ is disjoint from $B$ then $\delta^+(H) = \delta^-(H) = 0$,
because $A \cup C$ induces a transitive tournament.
On the other hand, in case $V(H) \cap B \neq \emptyset$, let $v$ be the vertex with zero indegree in the transitive tournament induced by $V(H) \cap B$.
Since \hbox{$\deg^-_H(v) \le \left|N^-_T(v) \cap C\right| = \ell$}, we conclude that
\[
\delta^\pm(H) \le  \ell = \frac{d+1}{2k} = \left(1+\frac{d+\ell}{dk}\right)\cdot \frac{d(d+1)}{2n} \le \left(1+\frac{k+1}{k^2}\right)\cdot \frac{d(d+1)}{2n}
\,.\]

It remains to consider the case $n > n_0$. 
Given the parameters $k$ and $\ell$, let $T_0$ be the previously constructed tournament on $n_0$ vertices.
In order to construct the desired tournament~$T$, add to $T_0$, one by one, a~total of~\mbox{$n-n_0$} vertices in such a way that every newly added vertex has the indegree zero.
The obtained tournament~$T$ satisfies $\delta^+(T) = \delta^+(T_0) = d$, and any $H \subseteq T$ with $\delta^-(H) > 0$ satisfies that $H \subseteq T_0$.
In~particular, it holds that $\delta^\pm(H) \le \ell$ for~every~$H \subseteq T$.
\end{proof}

\section{Dense digraphs}\label{sec:outro}
If $d=o(n)$, then Theorem~\ref{thm:main} and Proposition~\ref{prop:construction} determine,
up to a multiplicative factor $(1+o(1))$,
the best possible bound on $\max_{H\subseteq G} \delta^\pm(H)$ that holds for every $n$-vertex digraph $G$ with $\delta^+(G) \ge d$.
In the regime when $d = \alpha n$ for some fixed $\alpha > 0$, Theorem~\ref{thm:main} guarantees the existence of a~subgraph $H$ with $\frac{\delta^\pm(H)}n \ge \frac{\alpha^2}2$.

On the other hand, for any rational $\alpha \le \frac{1}{4}$,  let us fix $k := \lfloor \frac{1}{\alpha}\rfloor - 2$ and a~large enough integer $D$ such that
both $\ell := \frac{D}{2k}$ and $n := \frac D\alpha$ are integers.
Set $d:=2k\ell - 1 = D-1$ and note that for such a choice of the parameters, it holds that~\mbox{$(k+1)d+\ell < n$}.
The parameters $k$, $\ell$ and $n$ satisfy the assumptions of Proposition~\ref{prop:construction},
hence there exists an $n$-vertex tournament~$T$ with $\delta^+(T) = d = (\alpha-\frac1n)n$ and $\max\limits_{H\subseteq G} \delta^\pm(H) \le \ell$.
Therefore, for $\alpha \le \frac{1}{4}$, it holds that
\begin{equation}\label{eq:construction_dense}
 \max_{H\subseteq G} \frac{\delta^\pm(H)}n
\le \frac{\ell}n = \frac{\alpha}{2\left(\lfloor\frac1\alpha\rfloor-2\right)} \le \frac{\alpha}{2\left(\frac1{\alpha}-3\right)}
\le \frac{\alpha^2(1+12\alpha)}2
\,,
\end{equation}
which is very close to the bound proven in Theorem~\ref{thm:main} for $\alpha$ small.

\subsection{Very dense digraphs}\label{sec:densedigraphs}
\def\constDenseDG{0.7832}
We have observed that if $\alpha > 0$ is small, then the constant $\frac{\alpha^2}{2}$ in Theorem~\ref{thm:main} is almost optimal; see~\eqref{eq:construction_dense}.
In this section, we prove that if $\alpha \ge \constDenseDG$ then every $n$-vertex digraph $G$ with $\delta^+(G) = \alpha n$ has $H \subseteq G$ with $\frac{\delta^\pm(H)}{n} > \frac{\alpha^2}2$.

Let us start with introducing some additional notation.
Given an integer~$n$ and $\alpha \in (0,1)$, let $\cG_{n,\alpha}$ be the set of all $n$-vertex digraphs $G$ with $\delta^+(G) \ge \alpha n$.
Next, let
\[
h(n,\alpha) := \min\limits_{G \in \cG_{n,\alpha}} \max_{H \subseteq G} \frac{\delta^\pm(H)}n
\hspace{2em}\mbox{and}\hspace{2em}
h(\alpha) := \liminf\limits_{n\to\infty} \, h(n,\alpha)
 \,.\]
Note that $h(\alpha) \le \alpha$ by considering random digraphs of density~$\alpha$. 
Theorem~\ref{thm:main} yields that $h(\alpha) \ge \frac{\alpha^2}2$ by repeatedly removing vertices of~a~given graph that have small indegree or outdegree.
However, if $\alpha \ge \constDenseDG$ then one gets a~stronger lower bound by removing only those vertices that have a~small indegree.
More precisely, we assert that
\[h(\alpha) \ge 1 - \sqrt{3 - 4 \alpha + \alpha^2}\,,\]
which is larger than $\frac{\alpha^2}2$ for $\alpha > \sqrt{ 2\sqrt2+2}-\sqrt2 \approx \constDenseDG$.

To this end, fix an $n$-vertex digraph $G=(V,A)$.
Without loss of generality, every vertex $v \in V$ satisfies $\deg_G^+(v) = \lceil\alpha n\rceil$.
For the rest of this section, we will ignore floors and ceilings since the resulting $O(1/n)$ difference has no impact on the value of $h(\alpha)$.
Now remove one by one all its vertices that currently have the indegree smaller than $\beta n$ for~an~appropriately chosen $\beta < \alpha$.
This procedure terminates either when all the~vertices of $G$ have been removed, or, when the minimum indegree of the current subgraph is at least~$\beta n$.

Set $\beta := 1 - \sqrt{3 - 4 \alpha + \alpha^2}$.
Suppose the indegree-greedy procedure described in the previous paragraph is applied to $G$,
and let $\tau$ denote the proportion of the vertices of $G$ that has been removed so far.
Let $v_1,v_2,\ldots,v_{\tau n}$ be those vertices in the order of their removal, and,
for every $i\in[\tau n]$ we set $W_i:=\{v_1,v_2,\ldots,v_i\}$.
Double counting the number of arcs of $G$ yields that
\begin{align}\label{eq:greedydoublecount:digraph}
\begin{split}
\alpha
&=
\sum_{i=1}^{\tau n} \frac{|N_G^-(v_i) \cap W_i| + |N_G^-(v_i) \setminus W_i|}{n^2}
+
\sum_{v \in V\setminus W_{\tau n}} \frac{\deg_G^-(v)}{n^2}\\
&<
\frac{\tau^2}2 + \beta \tau + (1-\tau)
.
\end{split}
\end{align}
Let $f(\tau) := \frac{\tau^2}2 + \beta \tau + (1-\tau) - \alpha$.
A bit tedious calculation reveals that $f(\alpha-\beta)=0$, which together with \eqref{eq:greedydoublecount:digraph} implies that the procedure has to terminate after removing less than $(\alpha-\beta)n$ vertices from $G$.

Let $H\subseteq G$ be the remaining graph after the procedure terminates, and let $\tau_0$ be the total proportion of the removed vertices from $G$.
In other words, $H$ has $(1-\tau_0)n$ vertices.
It follows that every vertex of $H$ has the indegree at least $\beta n$, and the outdegree at least $(\alpha-\tau_0)n$.
Therefore,
\[
\frac{\delta^\pm(H)}n \ge \min\{ \alpha-\tau_0 , \beta \}  \ge \min\{ \alpha-(\alpha-\beta) , \beta \} = \beta
= 1 - \sqrt{3 - 4 \alpha + \alpha^2}
\,.
\]
In particular, it holds that $h(\alpha) \ge \max\left\{\frac{\alpha^2}2, 1 - \sqrt{3 - 4 \alpha + \alpha^2} \right\}$.
We believe it is an interesting open problem to find stronger estimates on~$h(\alpha)$.
\begin{problem}
For a given $\alpha \in (0,1)$, determine the value of $h(\alpha)$.
\end{problem}

\subsection{Orientations of very dense undirected graphs}\label{sec:denseorgraphs}
\def\constDenseOG{0.4528}

We say that a digraph $G$ is an \emph{oriented graph}, i.e., an orientation of the edge-set of a simple undirected graph,
if every pair of the vertices of $G$ is joined by at most one arc.
For any $n$-vertex oriented graph $G$ with $\delta^+(G) = \alpha n$ it must hold that $\alpha < 1/2$,
so at the first glance the previous subsection does not provide anything useful for $G$.
However, the double-counting argument \eqref{eq:greedydoublecount:digraph} can be easily improved
for oriented graphs since $\deg_G^-(v) < (1-\alpha)n$ for~every vertex $v \in V(G)$.
In the rest of this section, we show that this improvement yields the existence of a subgraph $H \subseteq G$ such that
\begin{equation}\label{eq:orgraph_bound}
\frac{\delta^\pm(H)}{n} \ge 1 - \alpha - \sqrt{3 - 8\alpha + 4\alpha^2}
\,.
\end{equation}
Observe that the right-hand side of~\eqref{eq:orgraph_bound} is larger than $\frac{\alpha^2}2$ for $\alpha > \alpha_1$,
where $\alpha_1$ is the only positive real root of $x^4+4x^3-16x^2+24x-8$.
We note that $\alpha_1 < \constDenseOG$.

Let us now introduce some additional notation, which will be analogous to the one in Section~\ref{sec:densedigraphs}.
For an integer $n$ and $\alpha \in (0,\frac{1}{2})$,
let $\widehat{\cG}_{n,\alpha}$ be the~set of~all $n$-vertex oriented graphs $G$ such that $\delta^+(G) \ge \alpha n$,
and 
\[
\widehat h(n,\alpha) := \min\limits_{G \in \widehat{\cG}_{n,\alpha}} \max_{H \subseteq G} \frac{\delta^\pm(H)}n
\hspace{2em}\mbox{and}\hspace{2em}
\widehat h(\alpha) := \liminf\limits_{n\to\infty} \, \widehat h(n,\alpha)
\,.\]
Since $\widehat{\cG}_{n,\alpha} \subseteq {\cG}_{n,\alpha}$, it holds that $\widehat h(\alpha)\ge h(\alpha)$.
Moreover, $\widehat h(\alpha) \le \alpha$ by considering, for example,
random subsets of the arc-set of a balanced tournament.

We are now ready to explain the existence of $H \subseteq G$ satisfying~\eqref{eq:orgraph_bound}.
Fix an~$n$-vertex oriented graph $G \in \widehat{\cG}_{n,\alpha}$.
Without loss of generality, we will again assume that $\deg_G^+(v) = \lceil\alpha n\rceil$ for every $v \in V(G)$, and in the rest, we will ignore floors and ceilings.
Now consider removing one by one all the vertices that have the current indegree smaller than $\widehat \beta n$ for some fixed $\widehat \beta > 0$.
As we have already mentioned, every vertex $v \in V(G)$ satisfies $\deg_G^-(v) < (1-\alpha)n$.
Thus the double-counting argument analogous to~\eqref{eq:greedydoublecount:digraph} yields~the~following:
\begin{equation}\label{eq:greedydoublecount:orgraph}
\alpha
<
\frac{\tau^2}2 + \widehat \beta \tau + (1-\tau)(1-\alpha)
\,.
\end{equation}
Let $H$ be the remaining subgraph after the indegree-greedy procedure terminates, and set
\[
\widehat \beta := 1 - \alpha - \sqrt{3 - 8\alpha + 4\alpha^2}
\quad \mbox{and} \quad
\widehat f(\tau):= \frac{\tau^2}2 + \widehat \beta \tau + (1-\tau)(1-\alpha) - \alpha
\,.
\]
We claim that $H$ satisfies~\eqref{eq:orgraph_bound}.
Indeed, since $\widehat f(\alpha-\widehat \beta) = 0$, the inequality \eqref{eq:greedydoublecount:orgraph} implies that $H$ has $(1-\tau_0)n$ vertices for some $\tau_0 < \alpha-\widehat \beta$.
In particular,
\[
\frac{\delta^\pm(H)}n \ge \min\{ \alpha-\tau_0 , \widehat \beta \}  \ge \min\{ \alpha-(\alpha-\widehat \beta) , \widehat \beta \}
= 1 - \alpha - \sqrt{3 - 8\alpha + 4\alpha^2}
\,.
\]

We have shown that $\widehat h(\alpha) \ge \max\left\{\frac{\alpha^2}2, 1 - \alpha - \sqrt{3 - 8\alpha + 4\alpha^2} \right\}$.
As in the previous subsection, we find interesting to establish stronger estimates on~$\widehat{h}(\alpha)$.
\begin{problem}
For a given $\alpha \in \left(0,\frac{1}{2}\right)$, determine the value of $\widehat{h}(\alpha)$.
\end{problem}

\subsection{Normalized minimum semidegree}\label{sec:denserelative}

Throughout this subsection, fix a large $G \in \cG_{n,\alpha}$.
A variant of the question we consider in this paper is to maximize, over all subgraphs $H \subseteq G$, the~value of~$\delta^\pm(H)$ normalized by the number of vertices of $H$.
In other words, our aim is now to bound the expression
\begin{equation*}\label{eq:relativemax}
\mu(G) := \max_{\emptyset \neq H\subseteq G} \frac{\delta^\pm(H)}{|V(H)|}
\,.
\end{equation*}

Let $\ell := \lceil 1/\alpha\rceil$ and $b \in \mathbb{N}$.
Recall that the $b$-blowup of a directed cycle of length $\ell$, which we will denote by $C_\ell[b]$,
is the oriented graph with $b\ell$ vertices partitioned into $\ell$ independent sets $V_0,V_1,\ldots,V_{\ell-1}$ of size $b$ each
such that $uw$ is an arc for $u \in V_i$ and $w \in V_{j}$ if and only if $j-i \equiv 1$ modulo~$\ell$.
A moment of thought reveals that $\mu(C_\ell[b]) \le 1/\ell \le \alpha$.

On the other hand, we claim that $\mu(G) \ge \frac{\alpha}2$ for every $G \in\cG_{n,\alpha}$.
The celebrated Caccetta-H\"aggkvist conjecture~\cite{bib:CH78} states that $G$ contains a directed cycle  of length at most $\lceil 1/\alpha \rceil$.
In particular, a consequence of the conjecture would be $\mu(G) \ge \frac1{\lceil 1/\alpha \rceil}$.
A result of Chv\'atal and Szemer\'edi~\cite{bib:CS83} weakens the conclusion of the conjecture and establishes a directed cycle of length less than $\frac2\alpha$,
thus a~subgraph $H\subseteq G$ with $\frac{\delta^\pm(H)}{|V(H)|} \ge \frac\alpha2$ always exists.
Note that the multiplicative constant $2$ from~\cite{bib:CS83} was later improved by Shen~\cite{bib:S98}.

The subgraph $H\subseteq G$ found in the previous paragraph is only of a~constant order.
It is natural to ask whether one can also find $H' \subseteq G$ with $\frac{\delta^\pm(H')}{|V(H')|} = \Theta(\alpha)$  while $|V(H')|$ is ``large''.
The answer is yes via a relatively standard application of the regularity method for directed graphs:
Given $G \in \cG_{n,\alpha}$, partition its vertex-set using a degree-form of the regularity lemma (see, e.g., \cite{bib:AS04}),
and then apply the above-mentioned result of Chv\'atal and Szemer\'edi to the so-called cluster graph between the dense regular pairs.
It follows that $G$ contains $\Theta\left(n^\ell\right)$ cycles of length $\ell$ for some $2 \le \ell \le \frac2{\alpha-\eps}$,
hence a classical supersaturation result~\cite{bib:E64} yields that $G$ contains $C_\ell[b]$ with $b = \polylog(n)$.

\section*{Acknowledgements}
We would like to thank the referees for their valuable comments,
and Maria Aksenovich, Felix C. Clemen and Benny Sudakov for fruitful discussions on the topics related to this project.


\begin{thebibliography}{9}
  \def\burlalt#1#2{\href{#2}{#1}}
  \def\doi#1{\burlalt{\textsc{doi}:\detokenize{#1}}{https://dx.doi.org/#1}}
  \setlength{\parskip}{3pt}
  \setlength{\itemsep}{3pt plus 0.3ex}

\bibitem{bib:AS04}N. Alon, A. Shapira: \emph{Testing subgraphs in directed graphs}, J. Comput. System Sci \textbf{69} (2004), 354-382. \doi{10.1016/j.jcss.2004.04.008}
\bibitem{bib:B-JG}J. Bang-Jansen, G. Gutin: \emph{Digraphs: Theory, Algorithms and Applications}, Springer, London (2009). \doi{10.1007/978-1-84800-998-1}
\bibitem{bib:CH78} L. Caccetta, R. H\"aggkvist, \emph{On minimal digraphs with given girth}, Congr. Numer. \textbf{21} (1978), 181--187.
\bibitem{bib:CS83} V. Chv\'atal, E. Szemer\'edi: \emph{Short cycles in directed graphs}, J. Combin. Theory Ser. B \textbf{35} (1983), 323--327. \doi{10.1016/0095-8956(83)90059-X}
\bibitem{bib:E64} P.~Erd\H{o}s: \emph{On extremal problems of graphs and generalized graphs}, {\em Israel J. Math.} \textbf{2} (1964), 183--190.
\bibitem{bib:E6465} P.~Erd\H{o}s: \emph{On an extremal problem in graph theory}, Colloq. Math. \textbf{13} (1964/65), 251--254. \doi{10.4064/cm-13-2-251-254}
\bibitem{bib:HMSSY13} H.~Huang, J.~Ma, A.~Shapira, B.~Sudakov, R.~Yuster: \emph{Large feedback arc sets, high minimum degree subgraphs, and long cycles in Eulerian digraphs}, Combin. Probab. Comput. \textbf{22} (2013), 859--873. \doi{10.1017/S0963548313000394}
\bibitem{bib:S98} J. Shen: \emph{Directed triangles in digraphs}, J. Combin. Theory Ser. B \textbf{74} (1998), 405--407. \doi{10.1006/jctb.1998.1839}

\end{thebibliography}
\end{document}